\documentclass[12pt]{article}
\usepackage{amsthm}
\usepackage{amsfonts}
\usepackage{mathrsfs}
\usepackage{amsmath}
\usepackage[numbers,sort&compress]{natbib}
\usepackage[xetex]{graphicx}
\newtheorem{thm}{Theorem}[section]

\newtheorem{rem}{Remark}[section]

\newtheorem{defn}{Definition}[section]

\newtheorem{ass}{Assumption}[section]

\numberwithin{equation}{section}
\allowdisplaybreaks[4]
\begin{document}
\date{}
\pagestyle{plain}
\title{Continuous-time Markov decision processes under the risk-sensitive average cost criterion}
\author{Qingda Wei$^1$, \   Xian Chen$^2$\thanks{The corresponding author.} \thanks{Wei's  email:   weiqd@hqu.edu.cn;  Chen's email:  chenxian@amss.ac.cn}
  \ \\ $^1$School of Economics and Finance, \\
  Huaqiao University, Quanzhou, 362021, P.R. China\\
  $^2$School of  Mathematical Sciences, \\
Peking University, Beijing, 100871, P.R. China}
\date{}
\maketitle \underline{}
\begin{abstract}
This paper studies continuous-time Markov decision processes under the risk-sensitive average cost criterion. The state space is a finite set, the action space is a Borel space,  the cost and transition rates are bounded, and the risk-sensitivity coefficient can take arbitrary positive real numbers. Under the mild conditions, we develop a new approach to establish the existence of a solution to the risk-sensitive average cost optimality equation and obtain the existence of an optimal deterministic stationary  policy.


\end{abstract}
\setlength{\baselineskip}{0.25in}

\section{Introduction} \label{intro}
Continuous-time Markov decision processes (CTMDPs) have wide applications, such as the  queueing systems,   control of the epidemic,  telecommunication, population processes, inventory control; see, for instance, \cite{put,kit95,guo09}.  The expected average cost criterion is a commonly used optimality criterion in the theory of CTMDPs and has  been widely studied under different sets of optimality conditions; see, for instance, \cite{put,guo09,w} and the references therein. The random costs incurred during the finite time interval are evaluated by the mathematical expectation
in the  definition of the expected average cost criterion. In other words, the expected average cost criterion assumes that
the decision-makers are risk-neutral. However,  different decision-makers may have different risk preferences in the real-world applications. Hence, it is necessary for us  to consider the attitude of a decision-maker towards the risk  in the definition of the average cost criterion. As is well known, the  utility function is an important tool to characterize the risk preferences of the decision-makers. In particular, the exponential utility function is a commonly used utility function and
has been  applied to reflect the risk attitudes  of  the  decision-makers towards  the random costs incurred  in the MDPs; see, for instance,  \cite{ca,ca10,di,ja} for discrete-time MDPs and \cite{ghosh} for CTMDPs. The average optimality criterion in \cite{ca,ca10,di,ja,ghosh}  is called risk-sensitive average cost criterion because
 the risk preferences of the decision-makers are taken into consideration. To the best of our knowledge, \cite{ghosh} is the first work  to study the risk-sensitive average cost criterion for CTMDPs. The state space is a denumerable set, the cost rate
 function is nonnegative and bounded, the transition rates are bounded and satisfy the irreducibility condition and some Lyapunov-like inequality, and the  risk-sensitivity coefficient of the exponential utility function is positive and satisfies some additional relation in \cite{ghosh}.

In this paper  we further  study the risk-sensitive average cost criterion  in the class of all randomized Markov policies for CTMDPs.  The state space is a finite set and the action space is a Borel space. The  cost rate function is bounded and allowed to take both nonnegative and negative values. The transition rates are bounded  and the risk-sensitivity coefficient  is allowed to take arbitrary positive real numbers. Under the irreducibility condition and the continuity and compactness conditions, we employ a new approach  to  establish the existence of a solution to the risk-sensitive average cost optimality equation, from which the existence of optimal policies  is shown.
More precisely,   we first introduce an auxiliary risk-sensitive first passage optimization problem and  obtain  the properties of the optimal value function of the risk-sensitive first passage problem (see Theorem \ref{thm3.1}). Then using the Feynman-Kac formula and the results on the risk-sensitive first passage optimization problem,  we show that the pair of the optimal value functions of the risk-sensitive average cost criterion and  the risk-sensitive first passage  problem is  a solution to the risk-sensitive average cost optimality equation and that there exists an optimal deterministic stationary policy in the class of all randomized Markov policies  (see Theorem \ref{thm3.2}). As far as we can tell,
 the risk-sensitive first passage optimization problem for CTMDPs is discussed for the first time in this paper. Moreover,
 since we remove the nonnegativity of the cost rate function, the  Lyapunov-like inequality imposed on the transition rates
 and the additional relation required for the positive risk-sensitivity coefficient in \cite{ghosh}, the optimality conditions
 in this paper are  weaker than those in \cite{ghosh} except that the state space is a finite set.
Furthermore, we  deal with   the risk-sensitive average cost criterion in a more general class of policies than that in \cite{ghosh}  which investigates this criterion in the class of all deterministic stationary policies.

 The rest of this paper is organized as follows. In Section 2, we introduce the decision model and the risk-sensitive average cost criterion. In Section 3, we give the optimality conditions and the main results whose proofs are  presented in Section 4.

\section{The decision model} \label{sec2}
The decision model we are concerned with is composed of the following components
$$\{S, A, (A(i), i\in S), q(j|i,a), c(i,a)\},$$
where the state space $S$  is a finite set endowed with the discrete topology, the action space $A$ is a Borel space with the
Borel $\sigma$-algebra $\mathcal{B}(A)$, and $A(i)\in\mathcal{B}(A)$ is the set of all admissible actions in state $i\in S$.
Let $K:=\{(i,a)|i\in S, a\in A(i)\}$ be the set of all admissible state-action pairs. The real-valued transition rate  $q(j|i,a)$ satisfies the following properties: (i) For each fixed $i,j\in S$, $q(j|i,a)$ is measurable in $a\in A(i)$; (ii) $q(j|i,a)\geq0$ for all $(i,a)\in K$ and $j\neq i$; (iii) $\sum_{j\in S}q(j|i,a)=0$ for all $(i,a)\in K$.  The  real-valued cost rate function $c(i,a)$ is  measurable in $a\in A(i)$ for each $i\in S$.

A continuous-time Markov decision process evolves as follows. A decision-maker observes continuously the state of a dynamical system. When the system is in state $i\in S$, an action $a\in A(i)$ is chosen by the decision-maker
according to some decision rule  and  such an intervention has the following consequences: (i) a cost is incurred at the rate $c(i,a)$; (ii) the system remains in the state $i$ for a random time following the exponential distribution with the tail function given by $e^{q(i|i,a)t}$, and then jumps to a new state $j\neq i$ with the probability $-\frac{q(j|i,a)}{q(i|i,a)}$  (we make a convention that $\frac{0}{0}:=0$).

Let $S_{\infty}:=S\cup\{i_{\infty}\}$ with an isolated point $i_{\infty}\notin S$, $\mathbb{R}_+:=(0,+\infty)$, $\Omega^0:=(S\times \mathbb{R}_+)^{\infty}$,
$\Omega:=\Omega^0\cup \{(i_0, \theta_1,i_1, \ldots, \theta_{m-1},i_{m-1}, \infty, i_{\infty},\infty, i_{\infty},\ldots)| i_0\in S, \ i_l\in S, \ \theta_l\in\mathbb{R}_+ \  {\rm for \ each} \ 1\leq l\leq m-1, \ m\geq2\}$, and $\mathcal{F}$  be the Borel $\sigma$-algebra of $\Omega$.  For each $\omega=(i_0, \theta_1, i_1, \ldots)\in \Omega$, define    $X_0(\omega):=i_0$, $T_0(\omega):=0$,  $X_m(\omega):=i_m$, $T_m(\omega):=\theta_1+\theta_2+\cdots+\theta_m$ for $m\geq1$, $T_{\infty}(\omega):=\lim_{m\to\infty}T_m(\omega)$, and the state process
$$\xi_t(\omega):=\sum_{m\geq0}I_{\{T_m\leq t<T_{m+1}\}}i_m+I_{\{t\geq T_{\infty}\}}i_{\infty}\ \ {\rm for} \ t\geq0,$$ where $I_D$ denotes  the indicator function of a set $D$. The process after $T_{\infty}$ is regarded to be absorbed in the state $i_{\infty}$. Hence, we write $q(i_{\infty}|i_{\infty}, a_{\infty})=0$, $c(i_{\infty}, a_{\infty})=0$, $A(i_{\infty}):=\{a_{\infty}\}$,  $A_{\infty}:=A\cup\{a_{\infty}\}$, where $a_{\infty}$ is an isolated point. Let $\mathcal{F}_t:=\sigma(\{T_m\leq s, X_m=i\}: i\in S, s\leq t, m\geq0)$ for $t\geq0$, $\mathcal{F}_{s-}:=\bigvee_{0\leq t<s}\mathcal{F}_t$, and $\mathcal{P}:=\sigma(\{D\times \{0\}, D\in \mathcal{F}_0\}\cup \{D\times (s,\infty), D\in \mathcal{F}_{s-}, s>0\})$ which denotes the $\sigma$-algebra of predictable sets on $\Omega\times [0,\infty)$ related to $\{\mathcal{F}_t\}_{t\geq0}$.

Now we introduce the definition of a randomized Markov policy below.
\begin{defn}
{\rm A $\mathcal{P}$-measurable transition probability $\pi(\cdot|\omega,t)$ on $(A_{\infty}, \mathcal{B}(A_{\infty}))$,
 concentrated on $A(\xi_{t-}(\omega))$ is called a  randomized Markov policy if there exists a kernel $\varphi$ on
 $A_{\infty}$ given $S_{\infty}\times [0,\infty)$ such that $\pi(\cdot|\omega,t)=\varphi(\cdot|\xi_{t-}(\omega),t)$.
  A policy  $\pi$ is said to be deterministic stationary if there exists a  function $f$ on $S_{\infty}$
  satisfying $f(i)\in A(i)$ for all $i\in  S_{\infty}$ and $\pi(\cdot|\omega, t)=\delta_{f(\xi_{t-}(\omega))}(\cdot)$, where $\delta_x(\cdot)$ is the Dirac measure concentrated at the point $x$.
}
\end{defn}

The set of all randomized Markov policies and the set of all deterministic stationary policies are denoted by $\Pi$ and $F$, respectively.

For any initial state $i\in S$ and any $\pi\in\Pi$, Theorem 4.27 in \cite{kit95} gives the existence of a unique probability measure $P_i^{\pi}$ on $(\Omega,\mathcal{F})$. Moreover,  the expectation operator with respect to $P_i^{\pi}$ is denoted by $E_i^{\pi}$.

Fix an arbitrary risk-sensitivity coefficient $\lambda>0$ throughout this paper.
For any   $i\in S$ and $\pi\in\Pi$, the  risk-sensitive average cost criterion is defined by
\begin{align*}
  J(i,\pi)=\limsup_{T\to\infty}\frac{1}{\lambda T}\ln E_i^{\pi}\left[e^{\lambda\int_0^T\int_Ac(\xi_t,a)\pi(da|\xi_t,t)dt}\right].
\end{align*}
The corresponding optimal value function is given by
\begin{align*}
  J^*(i):=\inf_{\pi\in\Pi}J(i,\pi) \ \ {\rm for \ all} \ i\in S.
\end{align*}
\begin{defn}
{\rm A policy $\pi^*\in\Pi$ is said to be   optimal if $J(i,\pi^*)=J^*(i)$ for all $i\in S$.
}
\end{defn}

The  main  goals of this paper are to give the conditions for the existence of optimal policies and to
develop  a new approach to  establish the existence of a solution to the risk-sensitive average cost optimality equation.

\section{The optimality conditions and main results}
In this section, we establish the existence of a solution to the risk-sensitive average cost optimality equation, from which the existence of optimal policies can be shown. To this end, we first introduce the following optimality conditions.
\begin{ass}\label{ass3.1}
\begin{itemize}
   \item[{\rm(i)}] For each $i\in S$,  the set $A(i)$ is compact.
   \item[{\rm (ii)}] For each $i,j\in S$, the functions $c(i,a)$ and $q(j|i,a)$ are continuous in $a\in A(i)$.
  \item[{\rm (iii)}] For each $f\in F$, the corresponding continuous-time  Markov chain $\{\xi_t,t\geq0\}$
  is irreducible,  which means that
  for any two states $i\neq j$, there exist different states $j_1=i$, $j_2,\ldots$, $j_m$ such that
  $q(j_2|j_1,f)\cdots q(j|j_m,f)>0$, where $q(j|i,f):=q(j|i,f(i))$.
\end{itemize}
\end{ass}
\begin{rem}
{\rm Assumptions \ref{ass3.1}(i) and \ref{ass3.1}(ii) are the standard continuity and compactness conditions which have been
widely used in CTMDPs; see, for instance, \cite{guo09,ghosh,w} and the references therein.
Moreover, Assumption \ref{ass3.1}(i) and the Tychonoff theorem imply that $F$  is compact and metrizable.
Assumption \ref{ass3.1}(iii) is the so-called irreducibility condition  which is commonly used in the average cost criterion; see, for instance, \cite{guo09} for the expected average case and \cite{ghosh} for the risk-sensitive average case.
}
\end{rem}
In order to prove the existence of optimal policies, we introduce the following notation.

For any fixed state $z\in S$, set $\tau_z:=\inf\{t\geq T_1:\xi_t=z\}$ with $\inf\emptyset:=\infty$. For each $i\in S$ and $f\in F$, let $c(i,f):=c(i,f(i))$. Below we introduce a risk-sensitive first passage optimization problem which has not been discussed in the existing literature.
For each $g\in \mathbb{R}:=(-\infty,\infty)$,  $i\in S$ and $f\in F$,   we define
\begin{eqnarray}\label{3-1}
h_{g}(i,f):=\frac{1}{\lambda}\ln E_i^f\left[e^{\lambda\int_0^{\tau_z}\left(c(\xi_t,f)-g\right)dt}\right] \ \ {\rm
 and} \ \  h_{g}^*(i):=\inf_{f\in F}h_{g}(i,f).
\end{eqnarray}
The function $h_{g}^*$ on $S$ is called the optimal value function of the risk-sensitive first passage problem.
Moreover, we set
\begin{align}\label{g-1}
 G:=\left\{g\in \mathbb{R}|h^*_{g}(z)\leq0\right\} \ \ {\rm and} \ \  \overline{g}:=\inf G.
\end{align}

Now  we state the first main result on the properties of the functions $h_g$ and $h^*_{g}$.
\begin{thm}\label{thm3.1}
Under  Assumption  \ref{ass3.1}, the following statements hold.
\begin{itemize}
  \item[{\rm (a)}] The set $G$ is nonempty.
  \item[{\rm (b)}] For each  $g\in \mathbb{R}$ and $f\in F$, the function $h_{g}(\cdot, f)$ on $S$ satisfies the following equations:
  \begin{eqnarray}\label{3-2}
\left\{\begin{array}{ll}
e^{\lambda h_{g}(i,f)}=Q(i,f,g)\left(q(z|i,f)+\sum_{j\in S\setminus\{i, z\}}e^{\lambda h_{g}(j,f)}q(j|i,f)\right) \\
e^{\lambda h_{g}(z,f)}=Q(z,f,g)  \sum_{j\in S\setminus\{z\}}e^{\lambda h_{g}(j,f)}q(j|z,f)
\end{array}\right.
\end{eqnarray}
for  all $i\in S\setminus\{z\}$, where we set $Q(i,f,g):=\int_0^{\infty}e^{\lambda (c(i,f)-g)s+q(i|i,f)s}ds$ and make a convention that  $0\cdot\infty:=0$.
\item[{\rm (c)}] For each  $g\in \mathbb{R}$ and $i\in S$, the function  $Q(i,a,g):=\int_0^{\infty}e^{\lambda (c(i,a)-g)s+q(i|i,a)s}ds$ is continuous in $a\in A(i)$. Moreover,
    $Q(i,a,g)\left(q(z|i,a)+\sum_{j\in S\setminus\{i, z\}}e^{\lambda h^*_{g}(j)}q(j|i,a)\right)$ (for $i\in S\setminus \{z\}$) and  $Q(z,a,g)\sum_{j\in S\setminus\{z\}}e^{\lambda h^*_{g}(j)}q(j|z,a)$ are lower semi-continuous in $a\in A(i)$ and $a\in A(z)$, respectively.
\item[{\rm (d)}] For each $g\in G$,  the function $h_{g}^*$ on $S$ satisfies the  following equations
 \begin{align}\label{eq}
\left\{\begin{array}{ll}
e^{\lambda h^*_{g}(i)}=\inf_{a\in A(i)}\left\{Q(i,a,g)\left(q(z|i,a)+\sum_{j\in S\setminus\{i, z\}}e^{\lambda h^*_{g}(j)}q(j|i,a)\right)\right\}\\
e^{\lambda h^*_{g}(z)}=\inf_{a\in A(z)}\left\{Q(z,a,g)\sum_{j\in S\setminus\{z\}}e^{\lambda h^*_{g}(j)}q(j|z,a)\right\}
\end{array}\right.
\end{align}
for all $i\in S\setminus\{z\}$. Moreover, there exists a policy $f_{g}\in F$ with $f_{g}(i)\in A(i)$ attaining the minimum of (\ref{eq}), and  for any $f_{g}\in F$ with $f_{g}(i)\in A(i)$ attaining the minimum of (\ref{eq}), we have
$h_{g}(i,f_{g})=h^*_{g}(i)\in \mathbb{R}$ and $Q(i,f_{g},g)<\infty$ for all $i\in S$.
\item[{\rm (e)}] We have $\overline{g} \in G$ and $h^*_{\overline{g}}(z)=0$.
\end{itemize}
\end{thm}
\begin{proof}
  See Section \ref{sec5}.
\end{proof}

Below we present the second main result on the risk-sensitive average cost optimality equation (\ref{op}) and the existence of optimal policies.
\begin{thm}\label{thm3.2}
  Suppose that Assumption  \ref{ass3.1} is satisfied. Let $\overline{g}$ and $h^*_{\overline{g}}$ be as in (\ref{3-1}) and (\ref{g-1}). Then we have
  \begin{itemize}
    \item[{\rm (a)}] The pair $(\overline{g}, h^*_{\overline{g}})\in \mathbb{R}\times B(S)$ satisfies the following equation:
    \begin{align}\label{op}
      \lambda \overline{g} e^{\lambda h^*_{\overline{g}}(i)}= \inf_{a\in A(i)}\left\{ \lambda c(i,a)e^{\lambda h^*_{\overline{g}}(i)}+\sum_{j\in S}e^{\lambda h^*_{\overline{g}}(j)}q(j|i,a)\right\}
    \end{align}
    for all $i\in S$, where $B(S)$ denotes the set of all real-valued functions on $S$.  Moreover,  there exists $f^*\in F$ with $f^*(i)\in A(i)$ attaining the minimum of (\ref{op}).
    \item[{\rm (b)}] For any  $f^*\in F$ with $f^*(i)\in A(i)$ attaining the minimum of (\ref{op}), we have $J^*(i)=J(i,f^*)=\overline{g}$ for all $i\in S$. Hence, the policy $f^*$ is  risk-sensitive average optimal.
  \end{itemize}
\end{thm}
\begin{proof}
  See Section \ref{sec5}.
\end{proof}
\begin{rem}\label{rem3.2}
{\rm (a)  In this paper we use a new approach to obtain the existence of a solution to the risk-sensitive average cost optimality equation (\ref{op}). Moreover, we  discuss  the risk-sensitive average cost criterion in the class of all randomized Markov policies whereas \cite{ghosh} restricts the study of this criterion to the class of all deterministic stationary policies.

(b) Theorem \ref{thm3.2}  establishes the existence of a solution to the risk-sensitive average cost optimality equation and the existence of optimal policies under the weaker conditions than those in  \cite{ghosh} except that the state space is a finite set in this paper. More precisely, we retain the irreducibility condition and the standard continuity and compactness conditions imposed in \cite{ghosh},
and remove the condition (A5) (i.e., the Lyapunov-like inequality)  in \cite{ghosh}. Moreover,
the   cost rate function $c$ is assumed to be nonnegative  and  bounded and the positive risk-sensitivity coefficient $\lambda$ is required to satisfy the relation that $\lambda \max_{(i,a)\in K}c(i,a)<b$ (for some constant $b>0$) in \cite{ghosh} whereas we allow  the cost rate function to take both nonnegative and negative values and there are no
restrictions  on the positive risk-sensitivity coefficient.

}
\end{rem}

\section{Proofs of Theorems \ref{thm3.1} and \ref{thm3.2}}\label{sec5}
In this section, we give the proofs of Theorems \ref{thm3.1} and \ref{thm3.2}.
\begin{proof}
[Proof of Theorem
 \rm \ref{thm3.1}]
(a)  Let  $M:=\max_{(i,a)\in K}c(i,a)$. Then we have $h_{M}(i,f)\leq 0$ for all $i\in S$ and $f\in F$, which implies $h^*_{M}(z)\leq0$. Hence,  the set $G$ is nonempty.

(b) Fix any $g\in \mathbb{R}$ and $f\in F$. By (\ref{3-1}), for any $i\in S\setminus\{z\}$,  we  obtain
\begin{align}\label{3-3}
 e^{\lambda h_{g}(i,f)}&=E_i^f\left[e^{\lambda\int_0^{\tau_z}\left(c(\xi_t,f)-g\right)dt}I_{\{\tau_z=T_1\}}\right]+
 E_i^f\left[e^{\lambda\int_0^{\tau_z}\left(c(\xi_t,f)-g\right)dt}I_{\{\tau_z>T_1\}}\right]\nonumber\\
 &=E_i^f\left[e^{\lambda\left(c(i,f)-g\right)T_1}I_{\{\tau_z=T_1\}}\right]+
 E_i^f\left[e^{\lambda\int_0^{T_1}(c(\xi_t,f)-g)dt}I_{\{\tau_z>T_1\}}E_i^f\left[e^{\lambda\int_{T_1}^{\tau_z}
 \left(c(\xi_t,f)-g\right)dt}\big|\xi_{T_1}\right]\right]\nonumber\\
 &=E_i^f\left[e^{\lambda\left(c(i,f)-g\right)T_1}I_{\{\tau_z=T_1\}}\right]+
 E_i^f\left[e^{\lambda (c(i,f)-g)T_1}I_{\{\tau_z>T_1\}}e^{\lambda h_{g}(\xi_{T_1},f)}\right]\nonumber\\
 &=\int_0^{\infty}e^{\lambda (c(i,f)-g)s}e^{q(i|i,f)s}ds\left(q(z|i,f)+\sum_{j\in S\setminus\{i,z\}}e^{\lambda h_{g}(j,f)}q(j|i,f)\right),
\end{align}
where the last equality is due to Proposition B.8 in \cite[p.205]{guo09}. On the other hand, using the similar arguments of (\ref{3-3}), we have
\begin{align*}
e^{\lambda h_{g}(z,f)}&=E_z^f\left[e^{\lambda (c(z,f)-g)T_1}I_{\{\tau_z>T_1\}}e^{\lambda h_{g}(\xi_{T_1},f)}\right]\\
&=\int_0^{\infty}e^{\lambda (c(z,f)-g)s}e^{q(z|z,f)s}ds\sum_{j\in S\setminus\{ z\}}e^{\lambda h_{g}(j,f)}q(j|z,f).
\end{align*}
Hence, part (b) follows from  the last equality and (\ref{3-3}).

(c)  Fix any $g\in \mathbb{R}$ and $i\in S$.  Let $\{a_n,n\geq1\}\subseteq A(i)$ be an arbitrary sequence converging to $a\in A(i)$. We deal with the cases $Q(i,a,g)<\infty$ and $Q(i,a,g)=\infty$  as follows.\\
Case 1:  $Q(i,a,g)<\infty$. Assumption \ref{ass3.1}(ii) gives $\lim_{n\to\infty}\lambda c(i,a_n)+q(i|i,a_n)=\lambda c(i,a)+q(i|i,a)$.  Note that $\lambda c(i,a)-\lambda g+q(i|i,a)<0$.   Thus, there exists a positive integer $n_0$ such that $\lambda c(i,a_n)-\lambda g+q(i|i,a_n)<0$ for all $n\geq n_0$. Hence, we obtain
$$ Q(i,a_n,g)=\frac{1}{\lambda g-\lambda c(i,a_n)-q(i|i,a_n)} \ \ {\rm for \ all} \ n\geq n_0,$$  which together with Assumption \ref{ass3.1}(ii) yields $\lim_{n\to\infty}Q(i,a_n,g)=Q(i,a,g)$. Therefore,
$Q(i,a,g)$ is continuous in $a\in A(i)$. \\
Case 2:   $Q(i,a,g)=\infty$.  The inequality $\limsup_{n\to\infty}Q(i,a_n,g)\leq Q(i,a,g)$ obviously holds. Thus, $Q(i,a,g)$ is upper semi-continuous in $a\in A(i)$. Moreover, by the Fatou lemma and Assumption \ref{ass3.1}(ii),  we have that $Q(i,a,g)$ is lower semi-continuous in $a\in A(i)$. Hence, $Q(i,a,g)$ is continuous in $a\in A(i)$.\\
Furthermore, it follows from Assumption \ref{ass3.1}(ii) and the Fatou lemma that $$Q(i,a,g)\bigg(q(z|i,a)+\sum_{j\in S\setminus\{i, z\}}e^{\lambda h^*_{g}(j)}q(j|i,a)\bigg) \ \ {\rm for}\ i\in S\setminus\{z\}$$ and  $Q(z,a,g)\sum_{j\in S\setminus\{z\}}e^{\lambda h^*_{g}(j)}q(j|z,a)$ are lower semi-continuous in $a\in A(i)$ and $a\in A(z)$, respectively.

(d) Fix any $g\in G$.  Employing (\ref{3-1}) and (\ref{3-2}), we get
\begin{eqnarray}\label{3-4}
\left\{\begin{array}{ll}
e^{\lambda h^*_{g}(i)}\geq\inf\limits_{a\in A(i)}\left\{Q(i,a,g)\left(q(z|i,a)+\sum_{j\in S\setminus\{i, z\}}e^{\lambda h^*_{g}(j)}q(j|i,a)\right)\right\}\\
e^{\lambda h^*_{g}(z)}\geq\inf\limits_{a\in A(z)}\left\{Q(z,a,g)\sum_{j\in S\setminus\{z\}}e^{\lambda h^*_{g}(j)}q(j|z,a)\right\}
\end{array}\right.
\end{eqnarray}
for all $i\in S\setminus\{z\}$.
Moreover,   by part (c) and Assumption \ref{ass3.1}(i),  there exists $f_{g}\in F$ with
$f_{g}(i)\in A(i)$ attaining the minimum of (\ref{3-4}) such that
\begin{eqnarray}\label{3-5}
\left\{\begin{array}{ll}
e^{\lambda h^*_{g}(i)}\geq Q(i,f_{g},g)\left(q(z|i,f_{g})+\sum_{j\in S\setminus\{i, z\}}e^{\lambda h^*_{g}(j)}q(j|i,f_{g})\right)\\
e^{\lambda h^*_{g}(z)}\geq  Q(z,f_{g},g)\sum_{j\in S\setminus\{z\}}e^{\lambda h^*_{g}(j)}q(j|z,f_{g})
\end{array}\right.
\end{eqnarray}
for all $i\in S\setminus\{z\}$.
For any $j\neq z$, Assumption \ref{ass3.1}(iii)  implies that there exist different states $j_1=z$, $j_2$, $\ldots$, $j_m=j$
such that $q(j_{n+1}|j_n,f_{g})>0$ for all $n=1,\ldots, m-1$, which together with $e^{\lambda h^*_{g}(z)}<\infty$ and (\ref{3-5}) yields $ e^{\lambda h^*_{g}(j)}<\infty$ for all $j\in S$. By
  (\ref{3-1}) and part (b) we obtain
\begin{eqnarray}\label{3-7}
\left\{\begin{array}{ll}
e^{\lambda h^*_{g}(i)}\leq Q(i,f_{g},g)\left(q(z|i,f_{g})+\sum_{j\in S\setminus\{i, z\}}e^{\lambda h_{g}(j,f_{g})}q(j|i,f_{g})\right) \\
e^{\lambda h^*_{g}(z)}\leq  Q(z,f_{g},g)\sum_{j\in S\setminus\{z\}}e^{\lambda h_{g}(j,f_{g})}q(j|z,f_{g})
\end{array}\right.
\end{eqnarray}
for all $i\in S\setminus\{z\}$.
On the other hand,  we have
\begin{align}\label{3-8}
  e^{\lambda h^*_{g}(i)}\geq &\sum_{m=1}^{n} E_i^{f_{g}}\left[e^{\lambda\int_0^{T_m}\left(c(\xi_t,f_{g})-g\right)dt}I_{\left\{\xi_{T_0}\neq z,\ldots,\xi_{T_{m-1}}\neq z, \xi_{T_m}=z\right\}}\right]\nonumber\\
  &+E_i^{f_{g}}\left[e^{\lambda\int_0^{T_n}\left(c(\xi_t,f_{g})-g\right)dt}e^{\lambda h^*_{g}(\xi_{T_n})}I_{\left\{\xi_{T_0}\neq z,\ldots, \xi_{T_n}\neq z\right\}}\right]
\end{align}
for all $i\in S\setminus\{z\}$ and $n=1,2,\ldots$. In fact, employing (\ref{3-5}), we obtain
\begin{align}\label{3-9}
   e^{\lambda h^*_{g}(\xi_{T_m})}\geq& E_i^{f_{g}}\left[e^{\lambda\int_{T_m}^{T_{m+1}}\left(c(\xi_t,f_{g})-g\right)dt}
   I_{\{\xi_{T_{m+1}}=z\}}\big|\xi_{T_m}\right]\nonumber\\
  &+E_i^{f_{g}}\left[e^{\lambda\int_{T_m}^{T_{m+1}}\left(c(\xi_t,f_{g})-g\right)dt}e^{\lambda h^*_{g}(\xi_{T_{m+1}})}I_{\{\xi_{T_{m+1}}\neq z\}}\big|\xi_{T_m}\right]
\end{align}
for all $\xi_{T_m}\in S\setminus\{z\}$ and $m=0,1,\ldots$. Thus, (\ref{3-8}) holds for $n=1$. Suppose that
(\ref{3-8}) holds for $n=l\geq1$. Then we have
\begin{align*}
  e^{\lambda h^*_{g}(i)}\geq & \sum_{m=1}^{l} E_i^{f_{g}}\left[e^{\lambda\int_0^{T_m}\left(c(\xi_t,f_{g})-g\right)dt}I_{\left\{\xi_{T_0}\neq z, \ldots,\xi_{T_{m-1}}\neq z, \xi_{T_m}=z\right\}}\right]\\
  &+E_i^{f_{g}}\left[e^{\lambda\int_0^{T_l}\left(c(\xi_t,f_{g})-g\right)dt}e^{\lambda h^*_{g}(\xi_{T_l})}I_{\left\{\xi_{T_0}\neq z, \ldots, \xi_{T_l}\neq z\right\}}\right]\\
  \geq & \sum_{m=1}^{l+1} E_i^{f_{g}}\left[e^{\lambda\int_0^{T_m}\left(c(\xi_t,f_{g})-g\right)dt}I_{\left\{\xi_{T_0}\neq z, \ldots,\xi_{T_{m-1}}\neq z, \xi_{T_m}=z\right\}}\right]\\
   &+E_i^{f_{g}}\left[e^{\lambda\int_{0}^{T_{l+1}}\left(c(\xi_t,f_{g})-g\right)dt}e^{\lambda h^*_{g}(\xi_{T_{l+1}})}I_{\left\{\xi_{T_0}\neq z, \ldots,  \xi_{T_{l+1}}\neq z\right\}}\right]
\end{align*}
for all $i\in S\setminus\{z\}$, where the last inequality is due to (\ref{3-9}). Hence, (\ref{3-8}) holds for $n=l+1$. Therefore, by the induction, we obtain that (\ref{3-8}) holds for all $n\geq1$. Moreover, employing (\ref{3-8}) we get
\begin{align}\label{s}
   e^{\lambda h^*_{g}(i)}\geq \sum_{m=1}^{\infty} E_i^{f_{g}}\left[e^{\lambda\int_0^{T_m}\left(c(\xi_t,f_{g})-g\right)dt}I_{\left\{\tau_z=T_m\right\}}\right]
   =e^{\lambda h_{g}(i,f_{g})},
\end{align}
 which together with (\ref{3-1}) implies
 \begin{align}\label{3-10}
   e^{\lambda h^*_{g}(i)}=e^{\lambda h_{g}(i,f_{g})}<\infty \ \ {\rm for \ all} \ i\in S\setminus\{z\}.
 \end{align}
 Thus, by  (\ref{3-7}) and (\ref{3-10}) we have
\begin{align*}
 e^{\lambda h^*_{g}(i)}\leq&Q(i,f_{g},g)\left(q(z|i,f_{g})+\sum_{j\in S\setminus\{i, z\}}e^{\lambda h^*_{g}(j)}q(j|i,f_{g})\right)\\
 =&\inf_{a\in A(i)}\left\{Q(i,a,g)\left(q(z|i,a)+\sum_{j\in S\setminus\{i, z\}}e^{\lambda h^*_{g}(j)}q(j|i,a)\right)\right\},
\end{align*}
which together with (\ref{3-4}) yields
\begin{align}\label{3-13}
   e^{\lambda h^*_{g}(i)}=\inf_{a\in A(i)}\left\{Q(i,a,g)\left(q(z|i,a)+\sum_{j\in S\setminus\{i, z\}}e^{\lambda h^*_{g}(j)}q(j|i,a)\right)\right\}
\end{align}
for all $i\in S\setminus\{z\}$.  Using the similar arguments of (\ref{3-10}) and (\ref{3-13}), we obtain
\begin{align}\label{c}
  e^{\lambda h^*_{g}(z)}=e^{\lambda h_{g}(z,f_{g})}<\infty \ \ {\rm and}\  \ e^{\lambda h^*_{g}(z)}= \inf_{a\in A(z)}\left\{Q(z,a,g)\sum_{j\in S\setminus\{z\}}e^{\lambda h^*_{g}(j)}q(j|z,a)\right\}.
\end{align}
Hence, the function $h^*_{g}$ on $S$ is a solution to the equation (\ref{eq}).  Furthermore, by (\ref{3-10}), (\ref{c}) and
Assumption \ref{ass3.1}(iii), we have $$e^{\lambda h^*_{g}(i)}=e^{\lambda h_{g}(i,f_{g})}\geq
E_i^{f_{g}}\left[e^{\lambda \left(\min_{(i,a)\in K}c(i,a)-g\right) \tau_z}\right]>0,$$ which implies $h^*_{g}(i)>-\infty$ for all $i\in S$.  Therefore,
from (\ref{3-10})-(\ref{c}), we conclude  that for any $f_{g}\in F$ with $f_{g}(i)\in A(i)$ attaining the minimum of (\ref{eq}), $h_{g}(i,f_{g})=h^*_{g}(i)\in \mathbb{R}$ and $Q(i,f_{g},g)<\infty$ for all $i\in S$.

(e)
Let $\{g_n,n\geq1\}\subseteq G$ be a  sequence satisfying
\begin{align}\label{3-14}
  g_n\geq g_{n+1} \ \ {\rm for   \ all} \ \ n\geq1 \ \ {\rm and} \ \ \lim_{n\to\infty}g_n=\overline{g}.
\end{align}
Then by  part (d), for each $n\geq1$, there exists $f_{g_n}\in F$ such that
\begin{eqnarray}\label{3-15}
\left\{\begin{array}{ll}
e^{\lambda h^*_{g_n}(i)}=Q(i,f_{g_n},g_n)\left(q(z|i,f_{g_n})
+\sum_{j\in S\setminus\{i, z\}}e^{\lambda h^*_{g_n}(j)}q(j|i,f_{g_n})\right)\\
e^{\lambda h^*_{g_n}(z)}=Q(z,f_{g_n},g_n)\sum_{j\in S\setminus\{z\}}e^{\lambda h^*_{g_n}(j)}q(j|z,f_{g_n})
\end{array}\right.
\end{eqnarray}
for all $i\in S\setminus\{z\}$. Since $F$ is compact, there exist a subsequence of $\{f_{g_n},n\geq1\}$ (still denoted by the same subsequence) and some $\widehat{f}\in F$ such that
\begin{align}\label{3-16}
f_{g_n}(i)\to \widehat{f}(i) \ \ {\rm as}  \ \ n\to\infty
\end{align}
for all $i\in S$. Moreover, using (\ref{3-1}) and (\ref{3-14}),  we have $h^*_{g_n}(i)\leq h^*_{g_{n+1}}(i)\leq h^*_{\overline{g}}(i)$  for all $n\geq1$,  which gives
\begin{align}\label{3-17}
  \lim_{n\to\infty}h^*_{g_n}(i)=:\widehat{h}(i)
  \leq h^*_{\overline{g}}(i) \ \ {\rm for \ all} \ i\in S.
\end{align}
Employing (\ref{3-14})-(\ref{3-17}) and the Fatou lemma, we obtain
\begin{eqnarray}\label{3-18}
\left\{\begin{array}{ll}
e^{\lambda \widehat{h}(i)}\geq Q(i,\widehat{f},\overline{g})\left(q(z|i,\widehat{f})
+\sum_{j\in S\setminus\{i, z\}}e^{\lambda \widehat{h}(j)}q(j|i,\widehat{f})\right)\\
e^{\lambda \widehat{h}(z)}\geq Q(z,\widehat{f},\overline{g})\sum_{j\in S\setminus\{z\}}e^{\lambda \widehat{h}(j)}q(j|z,\widehat{f})
\end{array}\right.
\end{eqnarray}
for all $i\in S\setminus\{z\}$.
Thus, by (\ref{3-18}) and the similar arguments of (\ref{s}), we get
 $\widehat{h}(i)\geq h_{\overline{g}}^*(i)$, which together with (\ref{3-17}) gives $\widehat{h}(i)=h_{\overline{g}}^*(i)$ for
all $i\in S$. Note that $\widehat{h}(z)\leq 0$. Hence, we have $h_{\overline{g}}^*(z)\leq 0$, which implies $\overline{g}\in G$.  Suppose that $h^*_{\overline{g}}(z)<0$. Let $f_{\overline{g}}\in F$ be the policy  with $f_{\overline{g}}(i)\in A(i)$ attaining the minimum of (\ref{eq}) and $\beta_n:=e^{n\lambda  h_{\overline{g}}^*(z)}$ ($n=1,2,\ldots$). By  part (d) we get  $\lambda c(i,f_{\overline{g}})-\lambda \overline{g}+q(i|i,f_{\overline{g}})<0$ for all $i\in S$. Thus, for each $n\geq1$, we define the new transition rates as follows:
\begin{align}\label{3-24}
  p_n(z|z,f_{\overline{g}}):=-\beta_{n+1}, \ p_n(j|z,f_{\overline{g}}):=-\frac{\beta_n e^{\lambda h^*_{\overline{g}}(j)}q(j|z,f_{\overline{g}})}{\lambda c(z,f_{\overline{g}})-\lambda \overline{g}+q(z|z,f_{\overline{g}})}\  {\rm for \ all} \ j\in S\setminus\{z\},
\end{align}
and for any $i\in S\setminus\{z\}$,
\begin{align}
  &p_n(i|i,f_{\overline{g}}):=-\beta_n e^{\lambda h^*_{\overline{g}}(i)}, \
  p_n(z|i,f_{\overline{g}}):=-\frac{\beta_n q(z|i,f_{\overline{g}})}{\lambda c(i,f_{\overline{g}})-\lambda \overline{g}+q(i|i,f_{\overline{g}})}, \label{3-25}\\
  &p_n(j|i,f_{\overline{g}}):=-\frac{\beta_n e^{\lambda h^*_{\overline{g}}(j)}q(j|i,f_{\overline{g}})}{\lambda c(i,f_{\overline{g}})-\lambda \overline{g}+q(i|i,f_{\overline{g}})} \ \ {\rm for \ all} \ j\in S\setminus \{i, z\}.\label{3-26}
\end{align}
For the policy $f_{\overline{g}}\in F$ and any initial state $i\in S$, the probability measure and expectation operator corresponding to the transition rates $p_n$ defined in (\ref{3-24})-(\ref{3-26}) are denoted by $P_{i,n}^{f_{\overline{g}}}$ and $E_{i,n}^{f_{\overline{g}}}$, respectively.  For any $\varepsilon>0$ and $n\geq1$,  define
$$H_{\varepsilon,n}(i):=\frac{1}{\lambda}\ln E_{i,n}^{f_{\overline{g}}}\left[e^{\lambda \varepsilon \tau_z}\right] \  {\rm for  \ all} \ i\in S.$$
By part (d), we have $e^{\lambda h^*_{\overline{g}}(i)}>0$ and $p_n(i|i,f_{\overline{g}})<0$ for all $i\in S$. Observe that $e^{\lambda h^*_{\overline{g}}(z)}<1$. Thus, there exists a positive integer $n_1$ such that
\begin{align}\label{3-28}
  \beta_{n_1}\leq \min_{i\in S}\left\{\left[\lambda \overline{g}-\lambda c(i,f_{\overline{g}})-q(i|i,f_{\overline{g}})\right]e^{-\lambda h^*_{\overline{g}}(i)}\right\}.
\end{align}
For  any  $\varepsilon\in \left(0,\min_{i\in S}\left\{-\frac{1}{\lambda}p_{n_1}(i|i,f_{\overline{g}})\right\}\right)=:O_{n_1}$, using (\ref{3-24})-(\ref{3-26}) and  the similar arguments of  part (b), we obtain
\begin{align}\label{3-27}
\left\{\begin{array}{ll}
  e^{\lambda H_{\varepsilon, n_1}(i)}=-\frac{1}{\beta_{n_1} e^{\lambda h^*_{\overline{g}}(i)}-\lambda \varepsilon}\left(\frac{\beta_{n_1} q(z|i,f_{\overline{g}})}{\lambda c(i,f_{\overline{g}})-\lambda \overline{g}+q(i|i,f_{\overline{g}})}+\sum\limits_{j\in S\setminus\{i, z\}}\frac{\beta_{n_1} e^{\lambda H_{\varepsilon,n_1}(j)+\lambda h^*_{\overline{g}}(j)}q(j|i,f_{\overline{g}})}{\lambda c(i,f_{\overline{g}})-\lambda \overline{g}+q(i|i,f_{\overline{g}})}\right)\\
  e^{\lambda H_{\varepsilon, n_1}(z)}=-\frac{1}{\beta_{n_1+1}-\lambda \varepsilon}\sum\limits_{j\in S\setminus\{z\}}\frac{\beta_{n_1}e^{\lambda H_{\varepsilon,n_1}(j)+\lambda h^*_{\overline{g}}(j)}q(j|z,f_{\overline{g}})}{\lambda c(z,f_{\overline{g}})-\lambda \overline{g}+q(z|z,f_{\overline{g}})}
\end{array}\right.
\end{align}
for all $i\in S\setminus\{z\}$.  On the other hand, by (\ref{3-24})-(\ref{3-26}) and Assumption \ref{ass3.1}, for each $i\in S$, we have
$P_{i,n_1}^{f_{\overline{g}}}(\tau_z<\infty)>0$. Set $\alpha_1:=\min_{i\in S}P_{i,n_1}^{f_{\overline{g}}}(\tau_z<\infty)$. Note  that $P_{i,n_1}^{f_{\overline{g}}}(\tau_z<\infty)=\lim_{n\to\infty}P_{i,n_1}^{f_{\overline{g}}}(\tau_z\leq n)$. Thus, for each  $i\in S$, there exists a positive integer $n(i)$ (depending on $i\in S$) such that
$P_{i,n_1}^{f_{\overline{g}}}(\tau_z\leq n(i))\geq P_{i,n_1}^{f_{\overline{g}}}(\tau_z<\infty)-\frac{\alpha_1}{2}\geq \frac{\alpha_1}{2}$. Hence, taking $t_1:=\max_{i\in S}n(i)$, we obtain
\begin{align*}
  P_{i,n_1}^{f_{\overline{g}}}(\tau_z>t_1)\leq 1-\frac{\alpha_1}{2} \ \ {\rm for \ all} \ i\in S.
\end{align*}
Employing the last inequality  and an induction argument, we get
\begin{align}\label{6-6}
  P_{i,n_1}^{f_{\overline{g}}}(\tau_z>nt_1)\leq \left(1-\frac{\alpha_1}{2}\right)^n
\end{align}
for all $i\in S$ and $n=1,2,\ldots$. Moreover, for any $\varepsilon_0\in O_{n_1}$ satisfying $\varepsilon_0<\frac{1}{\lambda t_1}\ln\frac{2}{2-\alpha_1}$, direct calculations  give
\begin{align}\label{w-5}
  e^{\lambda H_{\varepsilon_0,n_1}(i)}=&\sum_{m=0}^{\infty}E_{i,n_1}^{f_{\overline{g}}}\left[e^{\lambda \varepsilon_0 \tau_z}I_{\{\tau_z\in (mt_1,(m+1)t_1]\}}\right]\nonumber\\
  \leq&\sum_{m=0}^{\infty}e^{\lambda \varepsilon_0(m+1)t_1}E_{i,n_1}^{f_{\overline{g}}}\left[I_{\{\tau_z\in (mt_1,(m+1)t_1]\}}\right]\nonumber\\
  \leq& \sum_{m=0}^{\infty}e^{\lambda \varepsilon_0(m+1)t_1}P_{i,n_1}^{f_{\overline{g}}}(\tau_z>mt_1)\nonumber\\
  \leq&\sum_{m=0}^{\infty}e^{\lambda \varepsilon_0(m+1)t_1}\left(1-\frac{\alpha_1}{2}\right)^m\nonumber\\
  =&\frac{e^{\lambda\varepsilon_0 t_1}}{1-e^{\lambda\varepsilon_0t_1}(1-\frac{\alpha_1}{2})}<\infty
\end{align}
for all $i\in S$, where the third inequality follows from  (\ref{6-6}).
Choose  any  $\varepsilon_1\in(0,\varepsilon_0)$ satisfying $\varepsilon_1<\min\limits_{i\in S}\left\{\frac{1}{\lambda}\left[\lambda \overline{g}-\lambda c(i,f_{\overline{g}})-q(i|i,f_{\overline{g}})\right]\right\}$ and let
 $H_{\varepsilon_1, n_1}^*(i):=\beta_{n_1} e^{\lambda H_{\varepsilon_1,n_1}(i)+\lambda h^*_{\overline{g}}(i)}$ for all $i\in S$.
Then by  (\ref{3-28}) and  (\ref{3-27})  we  have
\begin{align*}
\left\{\begin{array}{ll}
 H^*_{\varepsilon_1, n_1}(i)\geq-\frac{1}{\lambda c(i,f_{\overline{g}})-\lambda \overline{g}+\lambda \varepsilon_1+q(i|i,f_{\overline{g}})}\left(\beta_{n_1} q(z|i,f_{\overline{g}})+\sum_{j\in S\setminus\{i, z\}} H^*_{\varepsilon_1, n_1}(j)q(j|i,f_{\overline{g}})\right)\\
  H^*_{\varepsilon_1, n_1}(z)\geq-\frac{1}{\lambda c(z,f_{\overline{g}})-\lambda \overline{g}+\lambda \varepsilon_1+q(z|z,f_{\overline{g}})}\sum_{j\in S\setminus\{z\}} H^*_{\varepsilon_1, n_1}(j)q(j|z,f_{\overline{g}})
\end{array}\right.
\end{align*}
for all $i\in S\setminus\{z\}$. By the last inequalities and the similar arguments of (\ref{s}),   we obtain
\begin{align}\label{3-30}
   H^*_{\varepsilon_1, n_1}(i)\geq\beta_{n_1}e^{\lambda h_{\overline{g}-\varepsilon_1}(i,f_{\overline{g}})}\geq \beta_{n_1}e^{\lambda h^*_{\overline{g}-\varepsilon_1}(i)}
\end{align}
for all $i\in S$. Let $\{\eta_m,m\geq1\}\subseteq(0,\varepsilon_1)$ be a sequence satisfying $\lim_{m\to\infty}\eta_m=0$.
By (\ref{w-5}) and the  dominated convergence theorem, we  have  $\lim_{m\to\infty}e^{\lambda H_{\eta_m,n_1}(z)}=1$. Thus, for any $\rho\in (0,e^{-\lambda h^*_{\overline{g}}(z)}-1)$, there exists a positive integer $m_0$ such that  $e^{\lambda H_{\eta_{m_0},n_1}(z)}<1+\rho$,  which implies $e^{\lambda H_{\eta_{m_0},n_1}(z)+\lambda h^*_{\overline{g}}(z)}<1$.
Moreover, it follows from (\ref{3-30})  that  $ h^*_{\overline{g}-\eta_{m_0}}(z)<0$. Hence, we obtain $\overline{g}-\eta_{m_0}\in G$, which leads to a contradiction that $\overline{g}\leq \overline{g}-\eta_{m_0}$. Therefore, we have $h^*_{\overline{g}}(z)=0$. This completes the proof of the theorem.
\end{proof}

Employing Theorem \ref{thm3.1} and the Feynman-Kac formula,  we prove Theorem \ref{thm3.2} below.

\begin{proof}
[Proof of Theorem
 \rm \ref{thm3.2}]
(a)   By   Theorems \ref{thm3.1}(d) and \ref{thm3.1}(e), we have that $(\overline{g}, h^*_{\overline{g}})\in \mathbb{R}\times B(S)$ satisfies the following equation
\begin{align}\label{w}
  e^{\lambda h^*_{\overline{g}}(i)}=\inf_{a\in A(i)}\left\{Q(i,a,\overline{g})\sum_{j\in S\setminus\{i\}}e^{\lambda h^*_{\overline{g}}(j)}q(j|i,a)\right\}
\end{align}
for all $i\in S$.  Moreover, it follows from the Weierstrass theorem in \cite[p.40]{al},
Theorem \ref{thm3.1}(c)  and Assumption \ref{ass3.1}(i) that there exists $f^*\in F$ with $f^*(i)\in A(i)$ attaining the minimum of (\ref{w}). Thus, we have
\begin{align}
  \lambda\overline{g}e^{\lambda h^*_{\overline{g}}(i)}=&\lambda c(i,f^*)e^{\lambda h^*_{\overline{g}}(i)}+\sum_{j\in S}
  e^{\lambda h^*_{\overline{g}}(j)}q(j|i,f^*)\label{w-1}\\
  \geq&\inf_{a\in A(i)}\left\{ \lambda c(i,a)e^{\lambda h^*_{\overline{g}}(i)}+\sum_{j\in S}e^{\lambda h^*_{\overline{g}}(j)}q(j|i,a)\right\}\label{w-2}
\end{align}
for all $i\in S$. Furthermore, employing (\ref{w}), we obtain
\begin{align}\label{w-3}
\lambda c(i,a)e^{\lambda h^*_{\overline{g}}(i)}+\sum_{j\in S}e^{\lambda h^*_{\overline{g}}(j)}q(j|i,a) \geq \lambda\overline{g}e^{\lambda h^*_{\overline{g}}(i)} \  \ {\rm for \ all} \ (i,a)\in K.
\end{align}
 In fact, if $\int_0^{\infty}e^{\left(\lambda c(i,a)-\lambda \overline{g}+q(i|i,a)\right)s}ds<\infty$, using (\ref{w}), we get
\begin{align*}
  \left(\int_0^{\infty}e^{\left(\lambda c(i,a)-\lambda \overline{g}+q(i|i,a)\right)s}ds\right)^{-1}e^{\lambda h^*_{\overline{g}}(i)} \leq\sum_{j\in S\setminus\{i\}}e^{\lambda h^*_{\overline{g}}(j)}q(j|i,a),
\end{align*}
which implies  (\ref{w-3}).
If $\int_0^{\infty}e^{\left(\lambda c(i,a)-\lambda \overline{g}+q(i|i,a)\right)s}ds=\infty$, we have $\lambda c(i,a)-\lambda \overline{g}+q(i|i,a)\geq0$. Then  we get
\begin{align*}
-\left(\lambda c(i,a)-\lambda \overline{g}+q(i|i,a)\right)e^{\lambda h^*_{\overline{g}}(i)}  \leq\sum_{j\in S\setminus\{i\}}e^{\lambda h^*_{\overline{g}}(j)}q(j|i,a),
\end{align*}
which gives (\ref{w-3}).  Hence, the assertion follows from (\ref{w-2}) and (\ref{w-3}).

(b) Fix any $f^*\in F$ with $f^*(i)\in A(i)$ attaining the minimum of (\ref{op}).
  By  the Feynman-Kac formula, we obtain
\begin{align*}
  &E_i^{f^*}\left[e^{\lambda\int_0^T\left(c(\xi_t,f^*)-\overline{g}\right)dt}e^{\lambda h^*_{\overline{g}}(\xi_T)}\right]-e^{\lambda h^*_{\overline{g}}(i)}\nonumber\\
  =&E_i^{f^*}\left[\int_0^Te^{\lambda \int_0^r \left(c(\xi_v,f^*)-\overline{g}\right)dv}\left(\left(\lambda c(\xi_r,f^*)-\lambda \overline{g}\right)e^{\lambda h^*_{\overline{g}}(\xi_r)}+\sum_{j\in S}
  e^{\lambda h^*_{\overline{g}}(j)}q(j|\xi_r,f^*)\right)dr\right],
\end{align*}
which together with (\ref{w-1}) yields
\begin{align*}
  E_i^{f^*}\left[e^{\lambda\int_0^T\left(c(\xi_t,f^*)-\overline{g}\right)dt}e^{\lambda h^*_{\overline{g}}(\xi_T)}\right]=e^{\lambda h^*_{\overline{g}}(i)}
\end{align*}
for all $i\in S$ and $T>0$.  Thus, using the last equality, we have
\begin{align*}
\frac{1}{\lambda T}\ln E_i^{f^*}\left[e^{\lambda\int_0^T\left(c(\xi_t,f^*)-\overline{g}\right)dt}\right]+\frac{1}{\lambda T}\ln \left(\min_{i\in S}e^{\lambda h^*_{\overline{g}}(i)}\right)-\frac{1}{T} h^*_{\overline{g}}(i)\leq \overline{g}
\end{align*}
for all $i\in S$ and $T>0$. Letting $T\to\infty$ in the last inequality, we obtain
\begin{align}\label{3-33}
 J^*(i)\leq J(i,f^*)\leq \overline{g}\ \ {\rm for \ all}  \ i\in S.
\end{align}
On the other hand,  for  any $\pi\in \Pi$ and $i\in S$,  the Feynman-Kac formula and (\ref{op}) yield
\begin{align*}
  &E_i^{\pi}\left[e^{\lambda\int_0^T\int_Ac(\xi_t,a)\pi(da|\xi_t,t)dt-\lambda \overline{g}T}e^{\lambda h^*_{\overline{g}}(\xi_T)}\right]-e^{\lambda h^*_{\overline{g}}(i)}\nonumber\\
  =&E_i^{\pi}\bigg[\int_0^Te^{\lambda\int_0^r\int_Ac(\xi_v,a)\pi(da|\xi_v,v)dv-\lambda \overline{g}r}\bigg(\left(\lambda \int_Ac(\xi_r,a)\pi(da|\xi_r,r)-\lambda \overline{g}\right)e^{\lambda h^*_{\overline{g}}(\xi_r)}\\
  &+\sum_{j\in S}e^{\lambda h^*_{\overline{g}}(j)}\int_Aq(j|\xi_r,a)\pi(da|\xi_r,r)\bigg)dr\bigg]\geq0
\end{align*}
for all $T>0$. Then employing the last inequality, we get
\begin{align}\label{3-36}
\overline{g}\leq  \frac{1}{\lambda T}\ln E_i^{\pi}\left[e^{\lambda\int_0^T\int_Ac(\xi_t,a)\pi(da|\xi_t,t)dt}\right]+\frac{1}{\lambda T} \ln\left(\max_{i\in S} e^{\lambda h^*_{\overline{g}}(i)}\right)-\frac{1}{T}h^*_{\overline{g}}(i)
\end{align}
for all $i\in S$, $\pi\in\Pi$ and $T>0$. Letting $T\to\infty$ in (\ref{3-36}), we have $\overline{g}\leq J(i,\pi)$ for all $\pi\in\Pi$, which gives
\begin{align}\label{3-37}
  \overline{g}\leq J^*(i) \ \ {\rm for \ all} \ i\in S.
\end{align}
Therefore, the desired result follows from (\ref{3-33}) and (\ref{3-37}).
\end{proof}

\end{document}